\documentclass[preprint,3p]{elsarticle}
\usepackage{amsmath, amsthm, amscd, amsfonts, amssymb, graphicx, color}

\newtheorem{cor}{Corollary}[section]

\newtheorem{thm}{Theorem}[section]
\newtheorem{lem}{Lemma}[section]
\newtheorem{rem}{Remark}[section]
\numberwithin{equation}{section}
\newproof{pf}{Proof}
\newproof{pot}{Proof of Theorem \ref{thm2}}
\newproof{poot}{Proof of Corollary \ref{co1}}
\numberwithin{equation}{section}
\newdefinition{ex}{Example}[section]
\journal{}
\begin{document}
\begin{frontmatter}

\title{Pointwise bounds for positive supersolutions of nonlinear elliptic problems involving the $p$-Laplacian and application}

\author{A. Aghajani $^{\text{a,b}}$ \corref{cor1}}
\cortext[cor1]{Corresponding author: A. Aghajani. Tel.
+9821-73913426. Fax +9821-77240472.} \ead{aghajani@iust.ac.ir}

\author{A.M. Tehrani $^{\text{a}}$\corref{}}
\ead{amtehrani@iust.ac.ir}

\address {$^{\text{a}}$School of Mathematics, Iran University of Science and
Technology, Narmak, Tehran 16844-13114, Iran.\\

$^{\text{b}}$School of Mathematics, Institute for Research in Fundamental Sciences (IPM), P.O.Box: 19395-5746, Tehran, Iran. }
\begin{abstract}
We  derive a priori bounds for positive supersolutions of  $ - \Delta_{p} u = \rho(x) f(u) $, where $p>1$ and $\Delta_{p}$ is the $p$-Laplace operator, in a smooth bounded domain of $R^{N}$ with zero
Dirichlet boundary conditions. We apply the results to  nonlinear elliptic eigenvalue problem $ - \Delta_{p} u = \lambda f(u) $, with Dirichlet boundary condition, where $ f $ is a nondecreasing continuous differentiable
function on $[0,\infty]$ such that $ f(0) > 0 $, $ f(t)^{\frac{1}{p-1}} $ is superlinear at infinity, and give sharp upper and lower bounds for the extremal parameter $ \lambda_{p}^{*} $. In particular, we consider the  nonlinearities $ f(u) = e^u $ and $ f(u) = (1+u)^{m} $ ($ m > p-1 $ ) and give explicit estimates on $ \lambda_{p}^{*} $. As a by-product of our results, we obtain a lower bound for the principal eigenvalue of the $ p $-Laplacian that improves obtained results in the recent literature for some range of $ p $ and $ N $.
\end{abstract}

\begin{keyword}
Nonlinear eigenvalue problem, Estimates of principal eigenvalue, Extremal parameter.
\end{keyword}

\end{frontmatter}
%
\section{Introduction}
%
Let $\Omega$ be a smooth bounded domain of $R^{N}$ and $p>1$. We consider the nonlinear elliptic problem
\begin{equation}\label{eq31}
\left\{\begin{array}{ll} -\Delta_{p} u =  \rho(x)f(u)  & {\rm }\ x \in \Omega,\\
\hskip7.3mm u \geqslant 0& {\rm }\ x \in \Omega,\\
\hskip7.3mm u = 0& {\rm }\ x \in \partial \Omega
\end{array}\right.
\end{equation}
where $\Delta_{p}$ is the $p$-Laplace operator defined by $\Delta_{p}u:=\text{div}(|\nabla u|^{p-2}\nabla u) $, $ \rho : \Omega \rightarrow \Bbb{R} $ is a nonnegative bounded function that is not identically zero and $ f $ satisfies \\
($ \mathcal{C} $) $ f : D_{f} = [0,a_{f}) \rightarrow \Bbb{R}^{+}:=[0,\infty) $ $ ( 0 < a_{f} \leqslant + \infty ) $ is a nondecreasing $C^{1}$ function with $ f(u) > 0 $ for $ u > 0 $.\\

We say that $u$ is a solution of (1.1) if $u\in W^{1,p}_{0}(\Omega)$, $u\geq0$ a.e., $\rho(x)f(u)\in L^{1}(\Omega)$, and
$$\int_{\Omega}|\nabla u|^{p-2}\nabla u.\nabla \varphi=\int_{\Omega}\rho(x)f(u)\varphi,~~\text{for~all}~~\varphi\in C^{\infty}_{c}(\Omega),$$
that is, for all $C^{\infty}$ functions $\varphi$ with compact support in $\Omega$. Note that, since $u$ is $p$-superharmonic we have that if $u\not \equiv0$ then $u>0$ a.e. in $\Omega$, by the strong maximum
principle (see \cite{CS, Mo99, Tr67, Va84}). A solution $u\in W^{1,p}_{0}(\Omega)$ is called a regular  solution of (1.1) if $\rho(x)f(u)\in L^{\infty}(\Omega)$. By the well known regularity results for degenerate
elliptic equations, if $u$ is a regular solution of (1.1) then $u\in C^{1,\alpha}(\bar{\Omega})$ for some $\alpha\in (0,1]$ (see for instance \cite{CS,Lie}). Also, we say that $u\in W^{1,p}_{0}(\Omega)$ is a supersolution of (2.1) if $\rho(x)f(u)\in L^{1}(\Omega)$ and
$-\Delta_{p} u \geq  \rho(x)f(u)$ in the weak sense. Reversing the inequality one defines the notion of subsolution.

The ball of radius $ R $ centered at $ x_{0} $ in $ \Bbb{R}^{N} $ will be denoted by $ B_{R} ( x_{0} ) $. Given a set $ \Omega $ in $ \Bbb{R}^{N} $ we let $ | \Omega | $ denote its $ N $-dimensional Lebesgue measure. The {\it $ p $-torsion} function $ \psi $ of a domain $ \Omega $ is the unique  solution of the problem
\begin{equation*}
\left\{\begin{array}{ll} -\Delta_{p} u = 1  & {\rm }\ x \in \Omega,\\
\hskip7.3mm u=0& {\rm }\ x \in \partial \Omega.
\end{array}\right.
\end{equation*}
We shall denote $\psi_M:= \sup_{x \in \Omega} \hskip1mm \psi (x)$.

In this paper, first we consider   $C^{1}$ positive supersolutions $u$ of $ (\ref{eq31}) $ in section 2  (by a positive solution we
mean a solution which is nonnegative and nontrivial) and give  explicit pointwise lower bounds for  $u$ under the condition that $ f $ satisfies ($\mathcal{C}$) and $f^{\frac{-1}{p-1}}\in L^{1}(0,a)$ for $a\in(0,a_{f})$. In particular, we prove that
\begin{equation*}
F(u(x)) \geqslant \frac{p-1}{p}~ \Big(\frac{\rho_{x}(d_{\Omega}(x))d_{\Omega}(x)^{p}}{N}\Big)^{\frac{1}{p-1}} ~~~ \textrm{for} ~ \textrm{all} ~ x \in \Omega,
\end{equation*}
where
$$ F(t) = \displaystyle{ \int_{0}^{t}\frac{ \textrm{d} s}{f(s)^{\frac{1}{p-1}}} }, ~ 0<t <a_{f},~~~\rho_{x}(r) = \inf \big\{ \rho(y) : ~ | y - x | < r \big\},~\text{and}~d_{\Omega}(x):=\text{dist}(x,\partial \Omega).$$
As an application, in section 3, we consider the eigenvalue problem
\begin{equation}\label{eq07}
\left\{\begin{array}{ll} - \Delta_{p}u = \lambda f(u) & {\rm }\ x \in \Omega, \\
\hskip7.3mm u = 0& {\rm }\ x \in \partial \Omega,
\end{array}\right.
\end{equation}
with $f$ satisfies ($ \mathcal{C} $) and define the extremal parameter $\lambda_{p}^{*}$ as
\begin{eqnarray*}
\lambda_{p}^{*}=\lambda_{p}^{*}(f,\Omega):=\sup \Big\{ \lambda>0 : \textrm{problem} ~ (\ref{eq07}) ~ \textrm{has} ~ \textrm{at} ~ \textrm{least} ~ \textrm{one} ~ \textrm{positive}~\textrm{bounded ~solution.} \Big\}.
\end{eqnarray*}

In the case when $f$, in addition, satisfies\\

(H) $ f : \Bbb{R}^{+} \rightarrow \Bbb{R}^{+} $ is $C^{1}$, $ f(0) > 0 $ and $ f(t)^{\frac{1}{p-1}} $ is superlinear at infinity (i.e., $ \lim_{ t \rightarrow \infty} \frac{g(t)}{t^{p-1}} = \infty $),\\

X. Cabr\'e and M. Sanch\'on in [\cite{CS}, Theorem 1.4] proved that  $\lambda_{p}^{*}\in (0,\infty)$ and for every $\lambda\in  (0,\lambda_{p}^{*})$  problem $ (\ref{eq07}) $ admits a minimal regular solution $u_{\lambda}$. Minimal means that it is smaller than any other supersolution of the problem. If in addition $f(t)^{\frac{1}{p-1}}$ is convex function satisfying $ \int_{0}^{\infty}\frac{ \textrm{d} s}{f(s)^{\frac{1}{p-1}}} < \infty $,
then $ (\ref{eq07}) $ admits no solution for $ \lambda > \lambda_{p}^{*}(f,\Omega)$. Moreover, the family $\{u_{\lambda}\}$ is increasing in $\lambda$ and every $u_{\lambda}$ is semi-stable in the sense that the second variation of the
energy functional associated with (1.2) is nonnegative definite (see Definition 1.1 in \cite{CS}). Using this property in \cite{CS} the authors established that $u^{*}=\lim_{\lambda \uparrow \lambda_{p}^{*}}u_{\lambda}$ is a solution of (1.2) with $\lambda=\lambda_{p}^{*}$ whenever $\liminf_{t\rightarrow\infty}\frac{tf'(t)}{f(t)}>p-1$; $u^{*}$ is
called the extremal solution.\\

 Let $ \lambda_{1} = \lambda ( p , \Omega ) $ be the first eigenvalue of $ p $-Laplacian subjected to Dirichlet boundary condition, i.e.,
\begin{equation}\label{eq24}
\lambda_{1}:= \min_{0 \neq v \in W^{1,p}_{0}(\Omega)}\frac{\int_{\Omega} | \nabla v |^{p} \textrm{d} x }{\int_{\Omega}| v|^{p} \textrm{d} x}.
\end{equation}
 Azorero and Peral in \cite{AP} showed that if $ f(u) = e^{u} $ then $ \lambda_{p}^{*} \leqslant \max \big\{ \lambda_{1},\lambda_{1} \big( \frac{p-1}{p} \big)^{p-1} \big\} $. Cabr\'e and Sanch\'on in \cite{CS} extended this result for every nonlinearity $ f $ satisfying (H), as
\begin{equation}\label{eq08}
\lambda_{p}^{*} \leqslant \max \Big\{ \lambda_{1}, \lambda_{1} \sup_{t \geqslant 0}\frac{t^{p-1}}{f(t)} \Big\}.
\end{equation}
In both proofs the authors (by a contradiction argument) used comparison principle  for the $ p $-Laplacian operator to construct, for every $\varepsilon>0$ sufficiently small, an increasing sequence of functions  whose limit is in $ W_{0}^{1,p}(\Omega) $ and solves the problem $ - \Delta_{p}w = ( \lambda_{1} + \varepsilon ) w^{p-1} $,  then used the fact that the first eigenvalue for the $ p $-Laplacian
is isolated to get a contradiction.\\
Before presenting our estimates on $\lambda_{p}^{*}$, first we improve $ (\ref{eq08}) $ as the following (using the homogeneity property of $ p $-Laplacian and  $ (\ref{eq08}) $ itself)
\begin{equation}\label{eq09}
\lambda_{p}^{*} \leqslant \lambda_{1} \sup_{t \geqslant 0}\frac{t^{p-1}}{f(t)}.
\end{equation}
Then we prove the following upper bound, without using the fact that  the first eigenvalue for the $ p $-Laplacian
is isolated,
\begin{equation*}
\lambda_{p}^{*} \leqslant \frac{1}{\psi_{M}^{p-1}} \Big( \int_{0}^{\infty} \frac{ \textrm{d} s}{f(s)^{\frac{1}{p-1}}} \Big)^{p-1},
\end{equation*}
where $\psi_{M}$ as defined before is the supremum (maximum) of the $p$-torsion function on $\Omega$. As we shall see, in many cases, this represents a sharper upper bound than $ (\ref{eq09}) $.

While there is no explicit formula for the lower bound in the literature for the critical parameter $ \lambda_{p}^{*} $ ($ p \neq 2 $), which is very important in application,  we shall prove the following lower bound for the extremal parameter of problem (1.2 ) with general nonlinearity $f$ satisfying $\mathcal{C}$, using the method of sub-super solution,
\begin{equation*}
\lambda_{p}^{*} \geqslant \max \Big\{\frac{1}{\psi_{M}^{p-1}}\sup_{0<t<a_{f}}\frac{t^{p-1}}{f(t)},~\sup_{0<\alpha<\frac{||F||_{\infty}}{\psi_M}} \alpha^{p-1}-\alpha^{p}\beta(\alpha) \Big\},
\end{equation*}
where
\begin{eqnarray*}
\beta(\alpha):=\sup_{x\in \Omega}f' \big( F^{-1}(\alpha \psi(x)) \big) f \big( F^{-1}(\alpha \psi(x)) \big)^{\frac{2-p}{p-1}} \big| \nabla \psi(x) \big|^{p},~~\text{and}~~||F||_{\infty}=\int_{0}^{a_{f}}\frac{ \textrm{d} s}{f(s)^{\frac{1}{p-1}}}
\end{eqnarray*}
In particular, if $ \Omega = B $ the unit ball in $ \Bbb{R}^{N} $ centered at the origin, then we have
\begin{equation}\label{eq10}
\lambda_{p}^{*} \geqslant \max \Big\{ N (\frac{p}{p-1})^{p-1}\sup_{0 < t < a_{f}}\frac{t^{p-1}}{f(t)}, ~ (\frac{p}{p-1})^{p-1} N \sup_{0 < \alpha < ||F||_{\infty}} \gamma_{\alpha} \Big\},
\end{equation}
where
\begin{eqnarray*}
\gamma(\alpha):=\alpha^{p-1} \Big( 1-\frac{p}{(p-1)N}\sup_{0 < t < a_{f}}f'(t)f(t)^{\frac{2-p}{p-1}}(\alpha-F(t)) \Big).
\end{eqnarray*}
As we shall see, the lower bound $ (\ref{eq10}) $, in some dimensions, gives the exact value of the extremal parameter for the standard nonlineareties $ f(u) = e^{u} $ and $ f(u) = (1+u)^{m} $ with ($ m > p-1 $). Moreover, when $ p = 2 $ the above bounds coincide with those given in \cite{AGT}. For example for the nonlinearity $ f(u) = e^{u} $ our results give
\begin{align*}
Np^{p-1} \geqslant \lambda_{p}^{*}(e^{u},B) \geqslant \left\{\begin{array}{ll} (\frac{p}{e})^{p-1}N& {\rm }\ N \leqslant \frac{p^{\frac{2p-1}{p-1}}}{e(p-1)},\\\\ (\frac{p-1}{p})^{p-1}\frac{N^{p}}{p}& {\rm }\ \frac{p^{\frac{2p-1}{p-1}}}{e(p-1)} < N \leqslant \frac{p^{2}}{p-1},\\\\  p^{p-1}(N-p) &{\rm }\ N > \frac{p^{2}}{p-1}. \end{array}\right.
\end{align*}
Also we show that our results can be used to estimate the first eigenvalue of $ p $-Laplacian from below. As it mentioned in \cite{KF}, while upper bounds for $ \lambda_{1} ( \Omega ) $ can be obtained by choosing particular test function $ v $ in $ (\ref{eq24}) $, but lower bounds are more challenging. For more details on estimates and asymptotics of the principal eigenvalue and eigenfunction of the $p$-Laplacian operator, we refer the reader to \cite{BD1,BD,BD2,KF}. For example when $ \Omega = B $ we shall prove the following lower bound, which is better than those given  in \cite{BD1,BD,KF}, for some range of $ p $ and $ N $ (see end of Section 3).

\begin{align*}
\lambda_{1}(B) \geqslant \left\{\begin{array}{ll} (\frac{p}{p-1})^{p-1}N& {\rm }\ N \leqslant \frac{p^{\frac{2p-1}{p-1}}}{e(p-1)},\\\\ (\frac{e}{p})^{p-1}\frac{N^{p}}{p}& {\rm }\ \frac{p^{\frac{2p-1}{p-1}}}{e(p-1)} < N \leqslant \frac{p^{2}}{p-1},\\\\  (\frac{pe}{p-1})^{p-1}(N-p) &{\rm }\ N > \frac{p^{2}}{p-1}. \end{array}\right.
\end{align*}
Finally in section 4, as an another application, we give a nonexistence result for positive supersolutions of  $ (\ref{eq31}) $ and apply this result to obtain upper bound for the pull-in voltage of a simple Micro-Electromechanical-Systems MEMS device.
\section{ Bounds for positive supersolutions of problem (1.1) }
In this section we consider positive supersolutions of problem $ (\ref{eq31}) $ and give pointwise lower bounds independent of any given supersolution under consideration. The following simple lemma is useful  in making bounds for solutions. The case $p=2$ is a variant of Kato's inequality used in \cite{BC,BCM}, see Lemma 1.7 in \cite{BC} and Lemma 2 in \cite{BCM}.
\begin{lem}\label{l1}
Let $ G : (0,a) \rightarrow \Bbb{R^{+}} $ ($a\leq\infty$) be an increasing concave $ C^{2} $ function and $u$ a continuously differentiable function on $\Omega$ with $0<u(x)<a$ for $x\in\Omega$. Then we have
\begin{equation*}
-\Delta_{p}G(u) \geq G'(u)^{p-1} (-\Delta_{p}u),~~x\in\Omega,
\end{equation*}
in the weak sense.
\end{lem}
\begin{pf}
For simplicity, we assume that $u$ is a $C^{2}$  function in $\Omega$. By smoothing $u$ and a standard argument one can prove it for a $C^{1}$ function $u$. Using the definition of $ \Delta_{p} $, the product rule for the divergence of  product of a scalar valued function and a vector field, $G'>0$ and $G''\leq0$ we simply compute
\begin{align*}
\Delta_{p}G(u) &=\textrm{div} \Big( \big| \nabla G(u) \big|^{p-2}\nabla G(u) \Big) \\
&=\textrm{div} \Big( G'(u)^{p-1} \big| \nabla u \big|^{p-2}\nabla u \Big) \\
&=\nabla \Big(G'(u)^{p-1} \Big) \centerdot \big| \nabla u \big|^{p-2} \nabla u + G'(u)^{p-1} \textrm{div} \Big( \big| \nabla u \big|^{p-2} \nabla u \Big) \\
&=(p-1)G''(u) G'(u)^{p-2} \nabla u \centerdot \big| \nabla u \big|^{p-2} \nabla u + G'(u)^{p-1} \Delta_{p}u \\
&=(p-1)G''(u) G'(u)^{p-2} \big| \nabla u \big|^{p}+G'(u)^{p-1} \Delta_{p}u\leq G'(u)^{p-1} \Delta_{p}u
\end{align*}
as desired. $ \square $
\end{pf}
Now let $ \psi_{\rho} $ be the unique solution of the equation
\begin{equation}\label{eq05}
\left\{\begin{array}{ll} -\Delta_{p}u=\rho(x) & {\rm }\ x \in \Omega, \\
\hskip7.3mm u = 0& {\rm }\ x \in \partial \Omega.
\end{array}\right.
\end{equation}
When $ \rho \equiv 1 $ then $ \psi_{1}=\psi $ is the \emph{$ p $-torsion} function of $ \Omega $ as in Section 1. Recall the definition of $\rho_{x}(r)$ as
\begin{eqnarray*}
\rho_{x}(r) := \inf_{ y \in B_{r}(x) } \rho(y) ~~~ 0 < r < d_{\Omega}(x)=\text{dist}(x,\partial\Omega).
\end{eqnarray*}
\begin{thm}\label{t1}
Let $ u $ be a $C^{1}$ positive supersolution of problem (1.1) with
$ f$ satisfies $\mathcal{C}$ and $ f^{\frac{1}{p-1}} \in L^{1}(0,a) $ for $ 0 < a < a_{f} $. Then
\begin{equation}\label{eq02}
F(u(x)) \geqslant \psi_{\rho}(x),  ~~~ x \in \Omega,
\end{equation}
where $F(0)=0$ and $ F(t) = \displaystyle{\int_{0}^{t} \frac{ \textrm{d} s}{f(s)^{\frac{1}{p-1}}}}, ~ t \in (0,a_{f}) $, and $\psi_{\rho}$ defined in $(2.1)$. Moreover, we have
\begin{equation}\label{eq03}
F(u(y)) \geqslant  \frac{p-1}{p}  \rho_{x} \big( d_{\Omega}(x) \big)^{\frac{1}{p-1}} ~ \frac{d_{\Omega}(x)^{\frac{p}{p-1}} - \big| x-y \big|^{\frac{p}{p-1}}}{N^{\frac{1}{p-1}}},  ~~~ |y-x| < d_{\Omega}(x).
\end{equation}
In particular,
\begin{equation}\label{eq04}
F(u(x)) \geqslant \frac{p-1}{p}~ \Big(\frac{\rho_{x}(d_{\Omega}(x))d_{\Omega}(x)^{p}}{N}\Big)^{\frac{1}{p-1}} ~~~ \textrm{for} ~ \textrm{all} ~ x \in \Omega.
\end{equation}
\end{thm}
\begin{pf}
First note that by the assumptions on $f$ and definition of $F$ we have $F'(t) = \displaystyle{\frac{1}{f(t)^{\frac{1}{p-1}}}} > 0 $ and $ F''(t) = \displaystyle{\frac{- f'(t)}{(p-1) f(t)^{\frac{p}{p-1}}}} \leqslant 0$, $0<t<a_{f}$, thus using Lemma $\ref{l1}$ (with $G=F$ and $a=a_{f}$) and the fact that $u$ is a supersolution, we can write
\begin{align*}
-\Delta_{p}F(u)&\geq F'(u)^{p-1}(-\Delta_{p}u) \\
&=\frac{1}{f(u)}(-\Delta_{p}u) \\
&\geq \rho(x)=-\Delta_{p}\psi_{\rho}.
\end{align*}
Now since we have $ F(u)=\psi_{\rho} = 0 $ on $ \partial \Omega $, then by the maximum principle we get $F \big( u(x) \big) \geqslant \psi_{\rho}(x)$ for every $x \in \Omega$ that proves $(\ref{eq02})$.\\
To prove $ (\ref{eq03}) $ we need to estimate $ \psi_{\rho} $ from below. Let $ x \in \Omega $. Then for $ y \in B_{d_{\Omega}(x)}(x) $ we get from $ (\ref{eq05}) $
\begin{equation}\label{eq11}
-\Delta_{p}\psi_{\rho}(y)=\rho(y) \geqslant \rho_{x} \big( d_{\Omega}(x) \big).
\end{equation}
Now consider the auxiliary function $w(y)=\big( \frac{p-1}{p} \big) \displaystyle{ \frac{d_{\Omega}(x)^{\frac{p}{p-1}}- \big| x-y \big|^{\frac{p}{p-1}}}{N^{\frac{1}{p-1}}}}$ which satisfies $-\Delta_{p} w=1$ in $B_{d_{\Omega}(x)}(x)$ and $w=0$ on $\partial B_{d_{\Omega}(x)}(x)$. Then from $ (\ref{eq11}) $ we get
\begin{eqnarray*}
-\Delta_{p}\psi_{\rho}(y) \geqslant -\Delta_{p} \Big( \rho_{x} \big( d_{\Omega}(x) \big)^{\frac{1}{p-1}}w(y) \Big),
\end{eqnarray*}
hence by the maximum principle $\psi_{\rho}(y) \geqslant \rho_{x} \big( d_{\Omega}(x) \big)^{\frac{1}{p-1}}w(y)$ in $B_{d_{\Omega}(x)}(x)$ that with the aid of $ (\ref{eq02}) $ proves $ (\ref{eq03}) $. Taking $y=x$ in $ (\ref{eq03}) $ gives $ (\ref{eq04}) $. $ \square $
\end{pf}
\section{Application to eigenvalue problem}
%
\subsection{Lower and upper bounds for $\lambda_{p}^{*}(f,\Omega)$}
Consider the nonlinear eigenvalue problem $ (\ref{eq07}) $. Before presenting our results based on Theorem $ \ref{t1} $, first we improve the upper bound $ (\ref{eq08}) $ for the extremal parameter $\lambda_{p}^{*}(f,\Omega)$ with $f$ satisfies (H), in the following lemma using the homogeneity property of $ p$- Laplacian and $ (\ref{eq08}) $ itself.
\begin{lem}
For the extremal parameter  of problem $ (\ref{eq07}) $ with $f$ satisfies (H), we have
\begin{equation}\label{eq14}
\lambda_{p}^{*} \leqslant \lambda_{1} \sup_{t \geqslant 0}\frac{t^{p-1}}{f(t)}.
\end{equation}
\end{lem}
\begin{pf}
 Assume that for some $\lambda>0$, $ u_{\lambda} $ be the minimal solution of $ (\ref{eq07}) $ and take an arbitrary positive number $ M \in ( 0 , \infty ) $. Then it is easy to see that the function $ w := M u_{\lambda} $ is a bounded solution of the equation
\begin{equation*}
\left\{\begin{array}{ll} -\Delta_{p} w = M^{p-1} \lambda g(w) & {\rm }\ x \in \Omega, \\
\hskip7.3mm w = 0& {\rm }\ x \in \partial \Omega,
\end{array}\right.
\end{equation*}
where $ g(u) := f ( \frac{u}{M} ) $. Hence from $ (\ref{eq08}) $ we must have
\begin{equation}\label{eq12}
M^{p-1} \lambda \leqslant \max \Big\{ \lambda_{1}, \lambda_{1} \sup_{t \geqslant 0}\frac{t^{p-1}}{g(t)} \Big\}.
\end{equation}
However, we have $\displaystyle{\sup_{t \geqslant 0}\frac{t^{p-1}}{g(t)}}=M^{p-1}\sup_{t \geqslant 0}\frac{t^{p-1}}{f(t)}$, thus from $ (\ref{eq12}) $ we get
\begin{equation}\label{eq13}
\lambda \leqslant \max \Big\{ \frac{\lambda_{1}}{M^{p-1}}, \lambda_{1} \sup_{t \geqslant 0}\frac{t^{p-1}}{f(t)} \Big\}.
\end{equation}
Now for $M$ sufficiently large we get from $ (\ref{eq13}) $ that
\begin{eqnarray*}
\lambda \leqslant \lambda_{1} \sup_{t \geqslant 0}\frac{t^{p-1}}{f(t)},
\end{eqnarray*}
which proves $ (\ref{eq14}) $. $ \square $
\end{pf}
\begin{thm}\label{t2}
Let $ \lambda_{p}^{*} $  be the extremal parameter  of problem $ (\ref{eq07}) $ with $f$ satisfy ($\mathcal{C}$). Then
\begin{equation}\label{eq15}
\lambda_{p}^{*}\leqslant \frac{1}{\psi_{M}^{p-1}} \Big( \int_{0}^{a_{f}}\frac{ \textrm{d} s}{f(s)^{\frac{1}{p-1}}} \Big)^{p-1},
\end{equation}
and
\begin{equation}\label{eq16}
\lambda_{p}^{*} \geqslant \max \Big\{\frac{1}{\psi_{M}^{p-1}}\sup_{0 < t < a_{f}}\frac{t^{p-1}}{f(t)}, ~ \sup_{0 < \alpha < \frac{ \| F \|_{\infty}}{\psi_\Omega}} \alpha^{p-1} - \alpha^{p} \beta ( \alpha ) \Big\},
\end{equation}
where $ \beta(\alpha):= \displaystyle{ \sup_{x \in \Omega}f' \Big( F^{-1} \big( \alpha \psi(x) \big) \Big) f \Big( F^{-1} \big( \alpha \psi(x) \big) \Big)^{\frac{2-p}{p-1}} \big| \nabla \psi(x) \big|^{p}} $.

In particular, if $ \Omega = B $ the unit ball in $ \Bbb{R}^{N} $, then we have
\begin{equation}\label{eq17}
\lambda_{p}^{*}  \geqslant \max \Big\{ N \big( \frac{p}{p-1} \big)^{p-1} \sup_{ 0 < t < a_{f} } \frac{t^{p-1}}{f(t)}, ~ \big( \frac{p}{p-1} \big)^{p-1} N \sup_{ 0 < \alpha < \| F \|_{\infty}} \gamma ( \alpha ) \Big\},
\end{equation}
where $ \gamma ( \alpha ) := \displaystyle{ \alpha^{p-1} \Big( 1 - \frac{p}{(p-1) N } \sup_{0 < s < F^{-1} ( \alpha ) }f'(s)f(s)^{\frac{2-p}{p-1}} \big( \alpha - F(s) \big) \Big) } $.
\end{thm}
\begin{pf}
From Theorem $ \ref{t1} $ (and, of course, with $ \rho \equiv 1 $ and $ f $ replaced by $ \lambda f $) we have $ F(u_{\lambda}(x)) \geqslant \lambda^{\frac{1}{p-1}} \psi(x) $, $ x \in \Omega $, thus
\begin{equation*}
\lambda^{\frac{1}{p-1}} \leqslant \frac{1}{\psi_M} \int_{0}^{u_{\lambda}(x_{0})} \frac{ \textrm{d} s}{f(t)^{\frac{1}{p-1}}} \leqslant \frac{1}{\psi_M} \int_{0}^{a_{f}} \frac{ \textrm{d} s}{f(t)^{\frac{1}{p-1}}},
\end{equation*}
that proves $ (\ref{eq15}) $.

We prove $ (\ref{eq16}) $ by the method of sub-supersolution. We construct a supersolution of $ (\ref{eq07}) $ in the form $ \bar{u} = \alpha \psi $ where $ \alpha > 0 $ is a scalar to be chosen later. We require that
\begin{eqnarray*}
\Delta_{p}\bar{u} + \lambda f(\bar{u}) = - \alpha^{p-1} + \lambda f ( \alpha \psi ) \leqslant 0, ~~~ in ~ \Omega.
\end{eqnarray*}
Since $ f $ is nondecreasing this is satisfied if $ \lambda \leqslant \frac{\alpha^{p-1}}{f(\alpha \psi_{M})} $ and making the optimal choice of $ \alpha $ we get the sufficient condition that $ \lambda \leqslant \frac{1}{\psi_{M}^{p-1}} \displaystyle{\sup_{0 < t < a_{f}}\frac{t^{p-1}}{f(t)} } $. On the other hand, $ \underline{u} = 0 $ is an allowable subsolution (note that we have $f(0)>0$), now  Proposition 2.1 in \cite{CS} implies that problem $ (\ref{eq07}) $ has a positive bounded solution, hence
\begin{equation}\label{eq18}
\lambda_{p}^{*} \geqslant \frac{1}{\psi_{M}^{p-1}} \sup_{0 < t < a_{f}}\frac{t^{p-1}}{f(t)}.
\end{equation}
Now we show that for $ \alpha \in (0,~\frac{||F||_{\infty}}{\psi_\Omega} ) $ the function $ \bar{\bar{u}}(x) = F^{-1}(\alpha \psi(x)) $ is a supersolution of $ (\ref{eq07}) $ for $ \lambda = \alpha^{p-1} - \alpha^{p} \beta ( \alpha ) $. To do this we simply compute $ \Delta_{p}\bar{\bar{u}}(x) $, using the facts that if we take $ y(t): = F^{-1} ( \alpha t ) $ then  $ \frac{dy}{dt} = \alpha f(y)^{\frac{1}{p-1}} $ and $ \frac{d^{2}y}{dt^{2}}= \frac{\alpha^{2}}{p-1}f'(y) f(y)^{\frac{3-p}{p-1}} $. We have
\begin{align*}
\Delta_{p} \bar{\bar{u}} (x) &= \Big( \alpha^{p} f' ( \bar{\bar{u}} ) f ( \bar{\bar{u}} )^{\frac{2-p}{p-1}} \big| \nabla \psi (x) \big|^{p} - \alpha^{p-1} \Big) f ( \bar{\bar{u}} ) \\
&\leqslant \Big( \alpha^{p} \sup_{x \in \Omega} f' ( \bar{\bar{u}} ) f ( \bar{\bar{u}} )^{ \frac{2-p}{p-1}} \big| \nabla \psi (x) \big|^{p} - \alpha^{p-1} \Big) f ( \bar{\bar{u}} ) \\
&= - \Big( \alpha^{p-1} - \alpha^{p} \beta(\alpha) \Big) f ( \bar{\bar{u}} ).
\end{align*}
In other words,
$ \Delta_{p} \bar{\bar{u}} (x) + \big( \alpha^{p-1}-\alpha^{p} \beta ( \alpha ) \big) f ( \bar{\bar{u}} ) \leqslant 0 $, and since we have $ \bar{\bar{u}} (x) = 0, ~ x \in \partial \Omega $, this shows that $ \bar{\bar{u}} $ is a supersolution of $ (\ref{eq07}) $ for $ \lambda = \alpha^{p-1} - \alpha^{p} \beta ( \alpha ) $. Using again
the fact that  $ \underline{u} = 0 $ is an allowable subsolution and Proposition 2.1 in \cite{CS}, we infer that problem $ (\ref{eq07}) $ with $ \lambda = \alpha^{p-1} - \alpha^{p} \beta ( \alpha ) $ has a positive bounded solution, hence
\begin{eqnarray*}
\lambda_{p}^{*} \geqslant \alpha^{p-1}-\alpha^{p}\beta(\alpha).
\end{eqnarray*}
Taking the supremum over $ \alpha \in (0,~\frac{||F||_{\infty}}{\psi_\Omega} ) $ and combining it with $ (\ref{eq18}) $, we obtain $ (\ref{eq16}) $.\\
When $ \Omega = B $ the unit ball of $ \Bbb{R}^{N} $, then we have the explicit formula $ \psi(x) = ( \frac{p-1}{p})\frac{1}{N^{\frac{1}{p-1}}}(1-|x|^{\frac{p}{p-1}}) $, hence $ \psi_{M} = \frac{p-1}{p} N^{\frac{-1}{p-1}} $ and $ \big| \nabla \psi(x) \big|^{p} = N^{\frac{-p}{p-1}} \big| x \big|^{\frac{p}{p-1}} $. Taking $ s = F^{-1} ( \alpha \psi(x) ) $ and make the change $ \alpha \rightarrow \frac{pN^{\frac{1}{p-1}}}{p-1}\alpha $ in $ (\ref{eq16}) $ we arrive at $ (\ref{eq17}) $. $\square$
\end{pf}
Now we compare $ (\ref{eq14}) $ with the upper bound for $\lambda_{p}^{*}$ in Theorem $ \ref{t2} $. First note that from $ (\ref{eq14}) $ and $ (\ref{eq16}) $ we get
\begin{equation}\label{eq19}
\frac{1}{\psi_{M}^{p-1}} \leqslant \lambda_{1}.
\end{equation}
Also, since $ f $ is nondecreasing we have $ || F ||_{\infty}^{p-1} = \Big( \displaystyle{ \int_{0}^{a_{f}} \frac{ \textrm{d} s}{f(s)^{\frac{1}{p-1}}} } \Big)^{p-1} \geqslant \displaystyle{ \sup_{0 < t < a_{f}}\frac{t^{p-1}}{f(t)} }:=\alpha_{f,p}$. Thus generally $ (\ref{eq15}) $ is better than $ (\ref{eq14}) $ if $\alpha_{f,p}||F||_{\infty}^{p-1} < \lambda_{1}\psi_{M}^{p-1}$.
However, in high dimension $ (\ref{eq15}) $ is much better than $ (\ref{eq14}) $, as one can show by the known results that $\lambda_{1}\psi_{M}^{p-1}\rightarrow \infty$ when $N \rightarrow \infty.$ For example, from \cite{KF, LW} if  $\Omega$ is a ball $B_{R}$ of radius $R$ then $\lambda_{1}(B_{R}) \geqslant (\frac{N}{pR})^{p}$, and since $\psi_{M}(B_{R}) = R^{\frac{p}{p-1}}(\frac{p-1}{p})N^{\frac{-1}{p-1}}$, then we have
\begin{eqnarray*}
\lambda_{1}\psi_{M}^{p-1} \geqslant \frac{(p-1)^{p-1}}{p^{2p-1}}  N^{p-1} \rightarrow \infty ~~~ \textrm{as} ~~~ N \rightarrow \infty.
\end{eqnarray*}
Another way to illustrate the sharpness of our results, we consider the quasilinear elliptic problem
\begin{equation}\label{eq20}
\left\{\begin{array}{ll} -\Delta_{p}u=\lambda f(u^{q}) & {\rm }\ x \in \Omega, \\
\hskip5.5mm u = 0& {\rm }\ x\in \partial \Omega,
\end{array}\right.
\end{equation}
where $ f : \Bbb{R^{+}} \rightarrow \Bbb{R}^{+} $ satisfies $\mathcal{C}$. The next theorem shows that $ (\ref{eq15}) $ and $ (\ref{eq16}) $ become sharp when $ q \rightarrow \infty $. We omit the proof as it follows along the same lines as that in the proof of the similar result for the case $ p = 2 $ in recent joint work of the authors with N. Ghoussoub \cite{AGT}.
\begin{thm}
The extremal parameter $\lambda^{*}_{p}=\lambda^{*}_{p}(f,\Omega,q)$  of problem $ (\ref{eq20}) $ satisfies
\begin{eqnarray*}
\lim_{q \rightarrow \infty} \lambda^{*}_{p}  = \frac{1}{f(0)\psi_{M}^{\frac{1}{p-1}}}
\end{eqnarray*}
In particular, when $ f(0) = 1 $ and $ \Omega $ is the unit ball $ B $ then
\begin{eqnarray*}
\lim_{q \rightarrow \infty} \lambda^{*}_{p} = \big( \frac{p}{p-1} \big)^{p-1} N.
\end{eqnarray*}
\end{thm}
\begin{ex}
Consider problem $ (\ref{eq07}) $ with $ f(u) = e^{u} $ and $ \Omega = B $. Here, we have $ \displaystyle{\sup_{0 < t < \infty}\frac{t^{p-1}}{f(t)} = \frac{(p-1)^{p-1}}{e^{p-1}}} $ and $ \| F \|_{\infty} = p-1 $, thus from $ (\ref{eq15}) $ we get
\begin{equation*}
\lambda_{p}^{*} \leqslant N p^{p-1}.
\end{equation*}
Moreover, it is easy to see that the function
$ f'(t) f(t)^{\frac{2-p}{p-1}} \big( \alpha-F(t) \big) $ is decreasing, hence takes its maximum value at $ t = 0 $. Thus, $ \beta ( \alpha ) = \alpha^{p-1} - \frac{p}{(p-1) N } \alpha^{p} $. Now from $ (\ref{eq17}) $ we get
\begin{align*}
\lambda_{p}^{*}(e^{u},B) \geqslant \left\{\begin{array}{ll} (\frac{p}{e})^{p-1}N& {\rm }\ N \leqslant \frac{p^{\frac{2p-1}{p-1}}}{e(p-1)},\\\\ (\frac{p-1}{p})^{p-1}\frac{N^{p}}{p}& {\rm }\ \frac{p^{\frac{2p-1}{p-1}}}{e(p-1)} < N \leqslant \frac{p^{2}}{p-1},\\\\  p^{p-1}(N-p) &{\rm }\ N > \frac{p^{2}}{p-1}. \end{array}\right.
\end{align*}
\end{ex}
\begin{rem}
 Garcia-Azorero, Peral and Puel \cite{GP1, GP2} considered problem $ (\ref{eq07}) $ for $f(u)=e^{u}$ in a general bounded domain $\Omega$ and proved that if $N< p+ \frac{4p}{p-1}$ then the extremal solution $u^{*}$ is bounded. Also, if $N \geqslant p+ \frac{4p}{p-1}$ and $\Omega=B$ they showed that
 \begin{equation*}
u^{*}(x) = - p \ln |x| ~~~ \textrm{and} ~~~ \lambda_{p}^{*} = p^{p-1}(N-p),
\end{equation*}
Hence the extremal solution is unbounded in this range, implies that $ \lambda_{p}^{*} \geqslant p^{p-1}(N-p) $ in every dimension $ N $. So from $ (\ref{eq18}) $ we see that our formula gives the exact value of $ \lambda_{p}^{*} $ as a lower bound for $ \frac{p^{2}}{p-1} < N $  (without knowing the exact formula of $ u^{*} $), and better lower bound for $ N< p + \frac{4p}{p-1} $.
\end{rem}
\begin{ex}
Consider problem $ (\ref{eq07}) $ with $ f(u) = \big( 1+u \big)^{m} $, $ m > p-1 $ and $ \Omega = B $. Then from $ (\ref{eq15}) $ we get
\begin{equation*}
\lambda_{p}^{*} \leqslant \big( \frac{p}{p-1} \big)^{p-1} N \Big( \int_{0}^{\infty}(1+s)^{\frac{-m}{p-1}} \Big)^{p-1} = \Big( \frac{p}{m+1-p} \Big)^{p-1} N.
\end{equation*}
Also, here we have $\displaystyle{\sup_{0 < t < \infty}\frac{t^{p-1}}{f(t)} = \big( p-1 \big)^{p-1} \big( m+1-p \big)^{m+1-p}m^{-m}}$ and $||F||_{\infty}=\displaystyle{ \frac{p-1}{m+1-p} }$. Moreover, it is easy to see that the function
$f'(t)f(t)^{\frac{2-p}{p-1}} \big( \alpha-F(t) \big)$ is decreasing, hence takes the maximum at $ t = 0 $. So $\beta(\alpha)=\alpha^{p-1}-\frac{pm}{(p-1)N}\alpha^{p}$. Now from $ (\ref{eq17}) $ we get
\begin{equation}\label{eq21}
 \lambda_{p}^{*}\big(( 1+u )^{m},B\big) \geqslant \left\{\begin{array}{ll} Nm^{-m}p^{p-1}(m+1-p)^{m+1-p}& {\rm }\ N \leqslant \frac{p^{\frac{2p-1}{p-1}}}{p-1}(\frac{m+1-p}{m})^{\frac{m+1-p}{p-1}},\\\\ (\frac{p-1}{m})^{p-1}(\frac{N}{p})^{p} & {\rm }\ \frac{p^{\frac{2p-1}{p-1}}}{p-1}(\frac{m+1-p}{m})^{\frac{m+1-p}{p-1}} < N \leqslant \frac{mp^{2}}{(p-1)(m+1-p)},\\\\ (\frac{p}{m+1-p})^{p-1}~\frac{m(N-p)-N(p-1)}{m+1-p} &{\rm }\ N > \frac{mp^{2}}{(p-1)(m+1-p)}. \end{array}\right.
\end{equation}
\end{ex}
\begin{rem}
A. Ferrero in \cite{F} ~ ( also see  \cite{CS})  by introducing the exact formula of $u^{*}$, i.e., the radial function
\begin{eqnarray*}
u^{*}(x) = |x|^{-\frac{p}{m-p+1}}-1, ~ \textrm{corresponds} ~ \textrm{to} ~\tilde{\lambda} = (\frac{p}{m+1-p})^{p-1} ~ \frac{m(N-p)-N(p-1)}{m+1-p},
\end{eqnarray*}
proved that if $ N > p \frac{4p}{p-1} $ and $ m > m_{\sharp}$ see
\cite{F,CS}  for definition of $m_{\sharp} ) $ then
$\lambda_{p}^{*}=\tilde{\lambda}$. Hence from $ (\ref{eq21}) $ we
see that our formula for lower bound gives the exact value of
$\lambda_{p}^{*}$ when $\frac{mp^{2}}{(p-1)(m+1-p)} < N$, and
better bounds for all other cases.
\end{rem}
\begin{ex}
Consider problem $ (\ref{eq07}) $ with $ f(u) = \big( 1-u \big)^{-m} $, $ m > p-1 $ and $ \Omega = B $. Then from $ (\ref{eq15}) $ we get
\begin{equation*}
\lambda_{p}^{*}\leqslant \big( \frac{p}{p-1} \big)^{p-1} N \Big( \int_{0}^{1}(1-s)^{\frac{m}{p-1}} \Big)^{p-1} = \Big( \frac{p}{m+p-1} \Big)^{p-1} N.
\end{equation*}
Also, here we have $\displaystyle{\sup_{0 < t < 1}\frac{t^{p-1}}{f(t)} = \big( p-1 \big)^{p-1} \big( m+p-1 \big)^{1-m-p}m^{m}}$ and $||F||_{\infty}=\displaystyle{ \frac{p-1}{m+p-1} }$. Moreover, it is easy to see that the function
$f'(t)f(t)^{\frac{2-p}{p-1}} \big( \alpha-F(t) \big)$ is decreasing, hence takes the maximum at $t=0$. So $\beta(\alpha)=\alpha^{p-1}-\frac{pm}{(p-1)N}\alpha^{p}.$ Now from $ (\ref{eq17}) $ we get
\begin{equation*}
 \lambda_{p}^{*}\big(( 1-u )^{-m},B\big) \geqslant \left\{\begin{array}{ll} Nm^{m}p^{p-1}(m+p-1)^{1-m-p}& {\rm }\ N \leqslant \frac{p^{\frac{2p-1}{p-1}}}{p-1}(\frac{m}{m+p-1})^{\frac{m+p-1}{p-1}},\\\\ (\frac{p-1}{m})^{p-1}(\frac{N}{p})^{p} & {\rm }\ \frac{p^{\frac{2p-1}{p-1}}}{p-1}(\frac{m}{m+p-1})^{\frac{m+p-1}{p-1}} < N \leqslant \frac{mp^{2}}{(p-1)(m+p-1)},\\\\ (\frac{p}{m+p-1})^{p-1}~\frac{m(N-p)-N(p-1)}{m+p-1} &{\rm }\ N > \frac{mp^{2}}{(p-1)(m+p-1)}. \end{array}\right.
\end{equation*}
\end{ex}

In order to get more explicit formulas for $ \lambda_{p}^{*} $, here we give explicit upper and lower bounds for $\psi_{M}$. Let
\begin{eqnarray}\label{eq28}
r_{\Omega}:=\sup_{x \in \Omega} d_{\Omega} (x),
\end{eqnarray}
be the Chebyshev radius of $ \Omega $. Also, let $ d:=\frac{1}{2} diam (\Omega) $. Find  $x_{0},x_{1} \in \Omega$ such that $ B_{r_{\Omega}}(x_{0}) \subseteq \Omega \subseteq B_{d}(x_{1}) $. Then by comparing the $p$-torsion function $\psi$ of $\Omega$ with the $p$-torsions of $B_{r_{\Omega}}(x_{0})$ and $B_{d}(x_{1})$, i.e., functions
$$(\frac{p-1}{p})N^{\frac{-1}{p-1}}(r_{\Omega}^{\frac{p}{p-1}}-|x-x_{0}|^{\frac{p}{p-1}}) ~~\text{and}~~
(\frac{p-1}{p})N^{\frac{-1}{p-1}}(d^{\frac{p}{p-1}}-|x-x_{0}|^{\frac{p}{p-1}}),$$
respectively, we get
\begin{equation}\label{eq22}
\big( \frac{p-1}{p} \big) N^{\frac{-1}{p-1}} r_{\Omega}^{\frac{p}{p-1}} \leqslant \psi_{M} \leqslant \big( \frac{p-1}{p}) N^{\frac{-1}{p-1}} \big( \frac{diam (\Omega)}{2} \big)^{\frac{p}{p-1}}.
\end{equation}
Also, the following lower bound for $\psi_{M}$ form \cite{DG} is better than that in $ (\ref{eq22}) $ whenever $ r_{\Omega} $ is small with respect to the volume $ | \Omega | $ of $ \Omega $. Let $ \tau_{p}(\Omega) $ be the $ p $-torsional rigidity
\begin{eqnarray*}
\tau_{p} ( \Omega ) := \int_{\Omega} \psi(x) \textrm{d} x,
\end{eqnarray*}
then from (\cite{DG}, Theorem 5.1) we have
\begin{equation}\label{eq23}
\tau_{p}(\Omega) \geqslant \big( \frac{p-1}{2p-1} \big) \frac{ | \Omega |^{\frac{2p-1}{p-1}}}{P(\Omega)^{\frac{p}{p-1}}},
\end{equation}
where $ P ( \Omega ) $ is the perimeter of $ \Omega $. Now using  $ \tau_{p}(\Omega) \leqslant \psi_{M} | \Omega | $, then from $ (\ref{eq23}) $ we obtain
\begin{equation*}
\psi_{M} \geqslant \frac{p-1}{2p-1} (\frac{ | \Omega | }{P ( \Omega ) } )^{ \frac{p}{p-1} }.
\end{equation*}
Hence from Theorem $ \ref{t2} $ we get the following explicit bounds for $ \lambda_{p}^{*} $.
\begin{cor}
Let $ \lambda_{p}^{*}$  be the extremal parameter  of problem $ (\ref{eq07}) $ with $f$ satisfy $\mathcal{C}$. Then
\begin{equation*}
\big( \frac{p}{p-1} \big)^{p-1}\frac{2^{p}N}{diam (\Omega)^{p}}\sup_{0 < t < a_{f}}\frac{t^{p-1}}{f(t)} \leqslant \lambda_{p}^{*} \leqslant \theta_{p,\Omega} \Big( \int_{0}^{a_{f}} \frac{ \textrm{d} s}{f(s)^{\frac{1}{p-1}}} \Big)^{p-1},
\end{equation*}
where
\begin{eqnarray*}
\theta_{p,\Omega}:=\min \Big\{ \big( \frac{p}{p-1} \big)^{p-1}\frac{N}{r_{\Omega}^{p}},~ \big( \frac{2p-1}{p-1} \big)^{p-1} \big( \frac{P(\Omega)}{ | \Omega | } \big)^{p}\Big\}.
\end{eqnarray*}
\end{cor}
\subsection{lower bound for the first eigenvalue of the $ p $-Laplacian}
Here we show that how our results can be applied to estimate the first eigenvalue of $ p $-Laplacian from below. First we recall some results from the literature.
Let $ h ( \Omega ) $ be the Cheeger constant of $ \Omega $, i.e.,
\begin{eqnarray*}
h ( \Omega ) := \inf_{ \Omega } \frac{ | \partial D | }{ | D | },
\end{eqnarray*}
with $ D $ varying over all smooth domain of $ \Omega $ whose boundary $ \partial D $ does not touch $ \partial \Omega $ and with $ | \partial D | $ and $ | D | $ denoting ($ n-1 $)- and $ n $-dimensional measure of $ \partial D $ and $ D $, see \cite{KF}.
The following lower bound from \cite{LW} is the extension of the same result for $ p = 2 $ proved by Cheeger, see \cite{C}.
\begin{equation}\label{eq25}
\lambda_{1} ( \Omega ) \geqslant \big( \frac{ h ( \Omega ) }{p} \big)^{p}, ~~~ p \in (1,\infty).
\end{equation}
If $ \Omega $ is a ball we know that $ h ( \Omega ) = \frac{N}{R} $, (see \cite{KF}) hence from $ (\ref{eq25}) $ we have
\begin{equation}\label{eq26}
\lambda_{1} ( B_{R} ) \geqslant \big( \frac{N}{pR} \big)^{p}, ~~~ p \in (1,\infty).
\end{equation}
The lower bound $ (\ref{eq26}) $ becomes sharp when $ p \rightarrow 1 $, as it is shown by   V. Friedman and B. Kawhol in \cite{KF} that $ \lambda_{1} ( \Omega ) $ converges to the Cheeger constant $ h ( \Omega ) $ when $ p \rightarrow 1 $. However, it is not sharp when $p\rightarrow\infty$, as from \cite{JL} we know that
\begin{equation*}
\lim_{ p \rightarrow \infty } ( \lambda_{1} )^{ \frac{1}{p}} ( \Omega ) = \frac{1} { r_{\Omega} },
\end{equation*}
where $ r_{\Omega} $ is defined in $ (\ref{eq28}) $. Hence, $ \lim_{ p \rightarrow \infty} ( \lambda_{1} )^{\frac{1}{p}} ( \Omega ) = \frac{1}{R} $, while the $ p $-th root of the right hand side of $ ( \ref{eq26} ) $ goes to zero when $ p \rightarrow \infty $.

Here, we give some lower bounds for $ \lambda_{1} $ using our results. First note that from $ (\ref{eq19}) $ and $ (\ref{eq22}) $ we have
\begin{equation}\label{eq34}
\lambda_{1} ( \Omega ) \geqslant \frac{1}{\psi_{M}^{p-1}} \geqslant \big( \frac{p}{p-1} \big)^{p-1} \big( \frac{2}{diam ( \Omega ) } \big)^{p} N.
\end{equation}
In particular, in the special case when $ \Omega $ is the ball $ B_{R} $ then
\begin{equation}\label{eq27}
\lambda_{1} ( B_{R} ) \geqslant \frac{1}{\psi_{M}^{p-1}} \geqslant \big( \frac{p}{p-1} \big)^{p-1} \frac{N}{R^{p}},
\end{equation}
which is  recently obtained by J. Benedikt and P. Der\'abek in \cite{BD1}.

The lower bound $ (\ref{eq27}) $ is better than $ ( \ref{eq26} ) $ when $ N < \frac{p^{\frac{2p-1}{p-1}}}{p-1} $, and also becomes sharp in both critical cases $ p \rightarrow 1 $ and $ p \rightarrow \infty $.  Also, the following lower bound for $ \lambda_{1} $, which is a  consequence of Example 2.1 and $ (\ref{eq14}) $, gives better bound on $ \lambda_{1}(B) $, for more values of $ p $ and $ N $.
\begin{align}\label{eq32}
\lambda_{1}(B) \geqslant \left\{\begin{array}{ll} (\frac{p}{p-1})^{p-1}N& {\rm }\ N \leqslant \frac{p^{\frac{2p-1}{p-1}}}{e(p-1)}, \\\\ (\frac{e}{p})^{p-1}\frac{N^{p}}{p}& {\rm }\ \frac{p^{\frac{2p-1}{p-1}}}{e(p-1)} < N \leqslant \frac{p^{2}}{p-1},\\\\  (\frac{pe}{p-1})^{p-1}(N-p) &{\rm }\ N > \frac{p^{2}}{p-1}. \end{array}\right.
\end{align}
Benedikt and Der\'abek in \cite{BD} also presented upper and lower bounds for  $\lambda_{1}(\Omega)$ on a bounded domain $\Omega\subset R^{N}$. In particular, when $\Omega=B$ they proved that
\begin{eqnarray}\label{eq33}
\lambda_{1}(B) \geqslant N p.
\end{eqnarray}
Comparing $ (\ref{eq32}) $ and $ (\ref{eq33}) $, one can easily check that when $ 1 < p \leqslant 2 $ the lower bound $ (\ref{eq32}) $ is better than $ (\ref{eq33}) $ in every dimension $N$. Also, when $ p > 2 $ the same is true when $ N \geqslant \dfrac{p^{\frac{p+1}{p-1}}}{e} $.
\section{ Nonexistence results}
Here we show that how one can apply Theorem $ \ref{t1} $ to prove nonexistence of positive solutions of differential inequalities involving $ p $-Laplacian.

Consider the differential inequality
\begin{equation}\label{eq29}
\left\{\begin{array}{ll} -\Delta_{p}u \geqslant \lambda \rho(x)f(u) & {\rm }\ x \in \Omega, \\
\hskip5.5mm u \geqslant 0& {\rm }\ x \in \Omega,\\
\hskip5.5mm u \in W^{1,p}_{0}(\Omega).
\end{array}\right.
\end{equation}
\begin{thm}
Let $f$ satisfy ($\mathcal{C}$), and $ \rho : \Omega \rightarrow \Bbb{R} $ is a nonnegative function that is not identically zero. Then
\begin{itemize}
\item[$ i) $] Inequality $ (\ref{eq29}) $ has no positive  $C^{1}$  solution if
\begin{equation}\label{eq30}
\lambda > \big( \frac{p}{p-1} \big)^{p-1} \frac{ N || F ||_{\infty}^{p-1}}{ \sup_{x \in \Omega} \big\{ \rho_{x} \big( d_{\Omega}(x) \big) d_{\Omega}(x)^{p} \big\} }.
\end{equation}
\item[$ ii) $] If $ \rho(x) = | x |^{\alpha} $, $ \alpha > 0 $ and $ \Omega = B_{R} $, then the same is true if
\begin{equation*}
\lambda > \Big( \frac{\alpha + p}{p-1} || F ||_{\infty} \Big)^{p-1} ( \alpha + N ) R^{ - ( \alpha + p ) }.
\end{equation*}
\end{itemize}
\end{thm}
\begin{pf}
\begin{itemize}
\item[i)] If $ (\ref{eq29}) $ has a  positive solution $ u $, then from $ (\ref{eq04}) $ in Theorem $ \ref{t1} $ (by replacing $ f $ with $ \lambda   f $) we get
\begin{eqnarray*}
\Big( \int_{0}^{u(x)} \frac{ \textrm{d} s}{f(s)^{\frac{1}{p-1}}} \Big)^{p-1} \geqslant \lambda \big( \frac{p-1}{p} \big)^{p-1} \rho_{x} \big( d_{\Omega}(x) \big) d_{\Omega}(x)^{p}, ~~~ x \in \Omega,
\end{eqnarray*}
and taking supremum on  both sides over $ \Omega $ we arrive at a contradiction with $ (\ref{eq30}) $.
\item[ii)]  Now, let $ \rho(x) = |x|^{\alpha} $ and $ \Omega = B_{R} $. In this case we can use $ (\ref{eq02}) $ directly. Indeed, it is easy to see that the function
\begin{eqnarray*}
\psi_{\rho}(x) = C \big( R^{\frac{\alpha+p}{p-1}} - | x |^{\frac{\alpha+p}{p-1}} \big), ~~~ \textrm{with} ~~~ C := \big( \frac{p-1}{\alpha+p} \big) \big( \alpha + N \big)^{\frac{-1}{p-1}},
\end{eqnarray*}
is the  solution of $ (\ref{eq05}) $ with  $ \rho(x) = | x |^{\alpha} $, hence from $ (\ref{eq02}) $ we must have
\begin{eqnarray*}
F(u(x)) \geqslant \lambda^{\frac{1}{p-1}}\psi_{\rho}(x),~x\in B_{R}.
\end{eqnarray*}
Taking supremum over $B_{R}$ we get the desired result. $\square$
\end{itemize}
\end{pf}
As an application of this result, consider the eigenvalue problem
\begin{equation*}
\left\{\begin{array}{ll} -\Delta u=\lambda \frac{|x|^{\alpha}}{(1-u)^{2}} & {\rm }\ x \in \Omega, \\
\hskip5.5mm u = 0& {\rm }\ x\in \partial \Omega,
\end{array}\right.
\end{equation*}
that in dimension $N=2$  models a simple \emph{Micro-Electromechanical-Systems} MEMS device, see \cite{CG,GH2,GH1,GPW}. Let $\lambda^{*}$ (called pull-in voltage) be the extremal parameter of the above eigenvalue problem, then from Theorem 2.7, we have
\begin{eqnarray*}
\lambda^{*} \leqslant \frac{(\alpha+2)(\alpha+N)}{3}R^{-(\alpha+2)}.
\end{eqnarray*}
This upper bound substantially improve the ones obtained in \cite{AGT,GH1,GPW}. It could be interesting to compare this bound to the lower bound for $\lambda^{*}$ given in \cite{GH2}, then we have
\begin{equation*}
\max \Big\{ \frac{4(\alpha+2)(\alpha+N)}{27},\frac{(\alpha+2)(3N+\alpha-4)}{9} \Big\} R^{-(\alpha+2)} \leqslant \lambda^{*} \leqslant \frac{(\alpha+2)(\alpha+N)}{3}R^{-(\alpha+2)}.
\end{equation*}
\section{Acknowledgement}
This research was in part supported by a grant from IPM (No.
94340123).

\end{document}